On skew polynomial rings

C.L.Wangneo

Jammu,J&K,India,180002

(E-mail:-wangneo.chaman@gmail.com )

Abstract

In this note we consider the links of prime ideals of certain skew polynomial rings and prove our main theorem, namely theorem [5], which states the following:- Let R be a noetherian ring that is link k-symmetric and let σ be an automorphism of R. Let S(R) denote the skew polynomial ring R[x,σ] . Let B be a prime ideal of S(R) that is extended from R, that is B=S(B ∩ R). Then , for a prime ideal D of S(R), there is a link D→B in the ring S(R) implies that D is an extended prime ideal of S(R).

Introduction

In this paper we study the links of prime ideals of certain skew polynomial rings S(R)=R[x,σ], where R is a noetherian ring and σ is an automorphism of R. Recall from [2] , page 178, or from [ 4] , there is a link from Q to P (Q→P for short) , for prime ideals Q and P of a noetherian ring R , if there is an ideal A of R such that QP≤A<Q∩P and the R- bi-module Q∩P/A is torsion free as a left R/Q - module and as a right R/P - module. To study the links of prime ideals of the skew polynomial ring R[x,σ] where R is a noetherian ring we essentially modify suitably the proof of the main theorem of [4] and to do so we make use of the definitions and the results of [ 5] . Then the main theorem that we prove states the following:-Let R be a noetherian ring that is link k-symmetric and let σ be an automorphism of R. Let S(R) denote the skew polynomial ring R[x,σ] . Let B be a prime ideal of S(R) that is extended from R, that is B=S(B ∩ R). Then , for a prime ideal D of S(R), there is a link D→B in the ring S(R) implies that D is an extended prime ideal of the ring S(R).

<u>Definitions and Notation:-</u>    We now mention the source of our reference regarding the definitions and notation that are in use throughout this paper. For the sake of convenience and relevance we will throughout this paper adhere to the definitions and notation as in   [4 ] and [ 5],wherein we had stated that all the definitions and notation of [4] and [5] had been adapted by us from [2] and [3]. Thus,for example, for an ideal I of a Noetherian ring R which is σ-semistable with respect to a fixed automorphism σ of R we denote by $I^o$  as in [3],lemma(6.9.9), the largest σ invariant ideal contained in I and which we also call, as in [5], the σ-invariant part of I. Clearly if I is a prime ideal of R then  $I^o$ is a σ- invariant semiprime ideal of R. For the definition of a skew polynomial ring ,see for example,[2],page 10. Moreover for all the basic results on skew polynomial rings that we shall need in this paper we again refer the reader  to [ 2] . Also  if  we denote by S(R) the  skew polynomial ring R[x,σ] , then for any ideal B of R we will denote by S(B) the set of all the polynomials in S(R) with coefficients from the ideal B. We will throughout assume the  known fact ,without further mention , that S(B) is not an ideal of the ring S(R) unless the ideal B is σ - invariant in which case S(B) equals the two sided ideal generated by the set B in the ring S(R). For a right module, M , over a ring we denote by |M| always the right krull dimension of the module M if this dimension exists. We must mention here that for a bimodule M over a ring R, we will again use the symbol |M| to denote  only the right krull dimension of M unless otherwise stated.  For the basic definition and results on krull dimension we

refer the reader to [1]. For a few more words about the terminology in this paper we mention that a ring R is noetherian means that R is a left as well as right noetherian ring. For two subsets A and B of a given set, A≤B means B contains A and A<B denotes A≤B but A≠B. For an ideal A of R ,c(A) denotes the set of elements of R that are regular modulo the ideal A. Also for a ring R,Spec.R denotes the set of prime ideals of R. Finally we mention that all our rings are with identity and all our modules are unitary.

Main theorem

We now consider the skew polynomial ring $S(R)=R[x,\sigma]$, where R is a noetherian ring and σ is an automorphism of R . We shall first prove results that describe certain properties of some prime ideals of the ring S(R). For this we start with the following definition from [4] .

Definition (1):- For a noetherian ring R we say a prime ideal P of R is link krull -symmetric (link k- symmetric for short ) if for any prime ideal Q of R such that Q is linked to P we have |R/P | =|R/Q|.

Proposition(2):- Let R be a noetherian ring and let let S(R) denote the skew polynomial ring R[x,σ]. Let B be a prime ideal in spec S(R) such that x is not an element of B. Then B ∩ R is a σ -stable (or σ-invariant) prime ideal of the ring R such that B ∩ R =$p^0$ for some σ- semistable prime ideal of R. In particular , if B is an extended prime ideal of S(R), then B=S($P^0$), for some σ- semistable prime ideal P of R.

Proof:- Use lemma(10.6.4) of [3 ] to observe that B is a σ- stable prime ideal of S(R).By lemma(10.6.5) of [3 ],we have that if R is a right noetherian σ- prime ring(using this definition and notation as in [3] chapter 10) then R is a semiprime ring. So $B \cap R$ is a σ-stable semiprime ideal of R. It is then not difficult to see that $B \cap R = P^0$ for some σ-semistable prime ideal P of R. To see this observe that if P is a minimal prime ideal over $B \cap R$ then $B \cap R \leq \sigma^i(p)$ for all $i \geq 1$. Moreover $\sigma^i(P)$ is also a minimal prime ideal over $B \cap R$. Since R is a right noetherian ring so R has only finitely many minimal prime ideals over $B \cap R$. As a consequence thus $\sigma^m(P) = P$ for some integer $m \geq 1$. Hence P is a σ-semistable prime ideal of R and $B \cap R \leq P^0$. We now show that $P^0 \leq B \cap R$. To see this assume $B \cap R \nleq P^0$, then as observed above since $B \cap R$ is a semiprime ideal so there is an ideal J of R such that $p0 \cap J = B \cap R$, where J is the intersection of all those minimal prime ideals over $B \cap R$ that do not contain $p0$. Since obviously then $(p0)(\sigma(J))$ is contained in the ideal $B \cap R$, so $\sigma(J)$ is contained in all those minimal prime ideals over $B \cap R$ that do not contain $p0$. Hence $\sigma(J)$ is contained in the ideal J. Clearly then we must have that J is also a σ-invariant ideal of R. Hence $S(P^0) \cap S(J) = S(B \cap R)$. Note however that

$S(B \cap R) \leq B$ where B is given to be a prime ideal of S(R). So $S(P^0) \cap S(J) \leq B$ and since $P^0$ is not contained in B by our assumption above, and also since B is a prime ideal of S(R), so we get that $S(J) \leq S(B \cap R)$. Thus $J \leq B \cap R$. Since as seen above $B \cap R \leq P^0$ so $J \leq P^0$. Thus $P^0 = B \cap R$ a contradiction to our assumption in the beginning. This proves the proposition. The rest of the proposition is obvious.

Remark:- In view of the above lemma, if B is merely an ideal of S(R) and x is not an element of B, then we cannot conclude that B is a σ-stable ideal of S(R).

<u>Proposition(3):-</u> Let R be a noetherian ring. Let σ be an automorphism of R and let S(R) denote the skew polynomial ring R[x,σ]. Let B be a prime ideal in spec S(R) which is extended from R, that is B=S(B ∩ R ). then the following hold true:-

1) B is a σ- stable prime ideal of S(R) and if there is a link D→B of prime ideals D,B of S(R) then D is also a σ- stable prime ideal of S(R).

2) Moreover in case of 1) above we can write B ∩ R = $P^0$ for some σ-semistable prime ideal P of R and D ∩ R = $Q^0$ for some σ- semistable prime ideal Q of R.

Proof:- It is clear that if B is an extended prime ideal of S(R), then x is not an element of B. Thus by proposition (2) above B is a σ- stable prime ideal of S(R). This proves 1). To prove 2) first observe that again by proposition (2) above we may write B ∩ R= $P^0$ for some σ-semistable prime ideal P of R. Next observe that since there is a link D→B of the prime ideals D and B of S(R) thus xƐB if and only if xƐD. But as seen above since x is not an element of B, hence x is not an element of D also. Thus again by proposition(2) above D is a σ- stable prime ideal of S(R). Clearly we may write ( as in the case of B) D ∩ R = $Q^0$ for some σ- semi stable prime ideal Q of R.

Now we prove theorem (4) below which essentially (after using proposition (3) proved above) characterises the existence of a link between two extended prime ideals of the skew polynomial ring S(R) where R is any general noetherian ring.

Theorem(4):- Let R be a noetherian ring. Let σ be an automorphism of R and let S ( R) denote the skew polynomial ring R[x,σ]. let Q and P be σ-semistable prime ideals of R. Let m≥1 be a common integer such that $σ^m(Q)=Q$ and $σ^m(P)=p$ and let $Q^0=∩σ^i(Q)$ and $P^0=∩σ^j(P)$. Then there is a link $σ^i(Q)→σ^j(p)$ of the prime ideals $\sigma^i(Q)$ and $\sigma^j(p)$ for some integers i,j($\geq$1) if and only if there is a link $S(Q^0) \to S(P^0)$ of the prime ideals $S(Q^0)$ and $S(p^0)$ of the ring S(R). Moreover a semistable prime ideal P of R is link k- symmetric if and only if the prime ideal $S(P^0)$ of the ring S(R) is link k-symmetric.

Proof:- To prove the theorem we adapt the proof of theorem (8) of [4]. Now let i,j be integers ($\geq$1) such that there is a link $\sigma^i(Q) \to \sigma^j(P)$ of the prime ideals $\sigma^i(Q)$ and $\sigma^j(P)$ of the ring R. Then by a similar arguement as that of theorem (8) of [5], this link induces a link $Q^0 \to P^0$ of the semiprime ideals $Q^0$ and $P^0$ of the ring R. As stated in the beginning, following the proof of [4], theorem (8), albeit some modification this gives a link $S(Q^0) \to S(P^0)$, of the prime ideals $S(Q^0)$ and $S(P^0)$ of the ring S(R). We now prove the converse. Suppose there is a link $S(Q^0) \to S(P^0)$ of the prime ideals $S(Q^0)$ and $S(P^0)$ of the ring S(R). Again following exactly the proof of [4], theorem (8), we get that there exists a link $Q^0 \to P^0$ of the semiprime ideals $Q^0$ and $P^0$ of the ring R. By [5], theorem (8), it follows that there exists a link $\sigma^i(Q) \to \sigma^j(P)$ of the prime ideals $\sigma^i(Q)$ and $\sigma^j(P)$ for some integers i,j($\geq$1) of the ring R. Moreover, if Q is linked to P and if P is link k-symmetric then $|R/\sigma^i(Q)|=|R/\sigma^j(P)|$, so using [3], theorem (6.5.4) we get that $|S(R)/S(Q^0)|=|S(R)/S(P^0)|$, proving that the ideal $S(P^0)$ is link k-symmetric. The proof of the converse is also similar.

The results proved so far now culminate in the proof of theorem (5) below which is our main theorem.

Theorem(5)(Main theorem):- let R be a noetherian ring let $\sigma$ be an automorphism of R. Let S(R) denote the skew polynomial ring R[x,$\sigma$]. Let B be a prime ideal of S(R) that is extended from R, that is, B=S(B $\cap$ R) and let B be link k-symmetric. Then, for a prime ideal D of S(R), there is a link D$\to$B in the ring S(R) implies that D is an extended prime ideal of S(R).

Proof:- We follow the proof of theorem (9) of [4]. Thus we suppose that

$D \cap B/L$ is the linking bimodule for the link $D \to B$ via an ideal L of S(R). Clearly since B is an extended prime ideal of S(R), so x is not an element of B. Thus as in proposition [2] above, we see that x is not an element of D also. Hence D and B are σ- stable prime ideals of S(R). Clearly x is not an element of B, and since $L \leq B$, so x is not an element of L also. Let $A = L \cap R$. By proposition (2) above we may write $B \cap R = P^0$ and $D \cap R = Q^0$ for some σ- semistable prime ideals P and Q of R. Now observe that $A \leq Q^0 \cap P^0$. Then as in theorem (9) of [4], we have two cases:-

Case(1):- $A < Q^0 \cap P^0$.

In this case we first observe by the same arguement as in [4], theorem (9) that the link $D \to B$ via the ideal L induces a link $D \cap R \to B \cap R$ in the ring R via the ideal $A = L \cap R$. Obviously we can consider this link as a link $Q^0 \to P^0$ of the semiprime ideals $Q^0$ and $P^0$ of the ring R via the ideal A. By [5], proposition (5), we choose A as the unique minimal ideal such that $A \leq L$ and such that $Q^0 \to P^0$ is a link via the ideal A. With this choice of A we now observe that A is a σ - invariant ideal of R. Let m (≥1) be a common integer such that $\sigma^m(Q) = Q$ and $\sigma^m(P) = P$. Then $A < Q^0 \cap P^0$ implies by an argument similar to [4], theorem (9), that there is a link $S(Q^0) \to S(P^0)$ of the prime ideals $S(Q^0)$ and $S(P^0)$ of the ring S(R) via the ideal S(A). Now again using an argument similar to [4], theorem (9), and using the hypothesis that $B = S(P^0)$ is link k-symmetric, we can show that $S(Q^0) = D$. Thus D is an extended prime ideal of the ring S(R).

case(2):- $A = Q^0 \cap P^0$. Three subcases arise in this case and we show that each subcase is an impossibility.

<u>Subcase(1):-</u> This case includes the cases $Q^0<P^0$ or $P^0<Q^0$. To see the impossibility of this subcase, we suppose without loss of generality that $P^0<Q^0$. Then by Goldie's theorem (see, for example [1]) we have that $c(P^0) \cap Q^0 \neq \phi$. Since $D \geq S(Q^0)$, and $B=S(P^0)$, so $D \cap c(B) \neq \phi$ which contradicts lemma (12.7) of [2], because we are given that D is linked to B. Hence $P^0<Q^0$ is an impossibility. Similarly we can show that $Q^0<P^0$ is an impossibility.

<u>Subcase (2) :-</u> $Q^0$ and $P^0$ are incomparable semi-prime ideals of R.

In this subcase first observe that A is a semiprime ideal and since the components of $Q^0$ and $P^0$ are finite sets of prime ideals so each such component is a minimal prime ideal over the ideal A. In particular the components Q and P of the semiprime ideal A are prime ideals minimal over the ideal A. Moreover if we assume that $m \geq 1$ is a common integer such that $\sigma^m(Q)=Q$ and $\sigma^m(P)=P$ then by the hypothesis of our subcase (2) we get that the sets of prime ideals

$\{Q, \sigma(Q), \sigma^2(Q), \ldots \sigma^{m-1}(Q)\}$ and $\{P, \sigma(P), \sigma^2(P), \ldots \sigma^{m-1}(P)\}$ are incomparable prime ideals that are minimal over the ideal A. Thus if we consider the ring $S(R)/S(A)$, then this is a semiprime ring with minimal prime ideal

$S(P^0)/S(A)=B/S(A)$ and since $S(A) \leq L$, so the link $D/S(A) \to B/S(A)$ via $L/S(A)$ is a link of the prime ideals $D/S(A)$ to the minimal prime ideal $B/S(A)$ of the semiprime ring $S(R)/S(A)$ which contradicts lemma (11.7) of [2]. Hence subcase (2) can not occur.

<u>Subcase(3):-</u> $Q^0=P^0$.

In this subcase $A = Q^0 \cap P^0 = P^0$. So $S(R)/S(A)$ is a semiprime ring such that $S(P)/S(A)$ is a minimal prime ideal of $S(R)/S(A)$. We use essentially the same argument as of subcase (2) above. To see this observe that the link $D \to B$ via the ideal $L$ induces a link $D/S(A) \to B/S(A)$ in the ring $S(R)/S(A)$ via the ideal $L/S(A)$. But in the prime ring $S(R)/S(A)$ since $B/S(A) = S(P)/S(A)$ is a minimal prime ideal of the ring $S(R)/S(A)$, then the existence of this link would again contradict lemma (11.7) of [2]. Hence we cannot have $Q^0 = P^0$.

Thus case(2) is impossible. Hence we conclude that $D = S(Q)$ as in case (1).

We now state an easy application of our main theorem to iterated skew polynomial rings in finitely many variables over a fully bounded noetherian ring R.

<u>Corollary(6):-</u> Let R be a fully bounded noetherian ring and let S(R) denote the skew polynomial ring $R[x,\sigma]$, in one variable x (where σ is an automorphism of R) over R. Let B be an extended prime ideal of S(R) and let cl.(B) denote the set of prime ideals of the ring S(R) that are linked to B. Then cl.(B) is a finite set and B is link krull symmetric.

<u>Proof:-</u> For the proof of this result first observe that if B is an extended prime ideal of S(R), so by proposition (3) above we can write $B = S(P^o)$ for some σ-semistable prime ideal P of R. Now observe that if D is a prime ideal of S(R) that is linked to B then by theorem(5) above D is also an extended prime ideal of S(R) and hence by proposition (3) above $D = S(P^o)$ for some σ-semistable prime ideal Q of R. Thus by [5], theorem (9), since Q and P are σ-semistable ideals, so D is linked to B if and only if some $\sigma^i(Q)$ is linked to the prime ideal P. Now define a map f from the set cl.P to cl.B, by defining $f(Q)=S(Q^o)$, where Q is a prime ideal linked to P. Now note first that by theorem (14.22) of [2], cl.P is a finite set, and hence Q must be a σ-semistable prime ideal of R. Now the prime ideal Q is linked to the prime ideal P implies by theorem (4) above that the prime ideal $S(Q^o)$ is linked to the prime ideal $B=S(P^o)$. Clearly thus f is a well defined map which obviously is a surjection. Since cl.P is a finite set, hence cl.(B) is also a finite set. Also since R is a fully bounded noetherian ring so the prime ideal P is link krull smmetric. Hence using [3], theorem (6.5.4) we get that B is also link krull symmetric. This concludes the proof of the corollary.

<u>Corollary (7):-</u> Let R be a fully bounded noetherian ring and let $S_n(R)$ denote the iterated skew polynomial ring in n-variables $x_i$; ( with respect to n-automorphisms $\sigma_i$ (1≤i≤n)) over the ring R. Let B be an extended prime ideal of $S_n(R)$. Let cl.(B) denote the set of all the prime ideals of the ring $S_n(R)$ that are linked to B, then cl.(B) is a finite set. Moreover B is link k-symmetric.

<u>Proof:-</u> The proof is by induction on n. Observe that if in the case n=0, we represent the skew polynomial ring $S_0(R)$ by the ring R then the result follows by theorem (14.22) of [2]. The case n=1 is the corollary (6) proved above. Now assume the induction hypothesis and observe that the ring $S_{n-1}(R)$ shares the following properties of the ring R, namely (1) for a prime ideal P of the ring $S_{n-1}(R)$ if cl. P is a finite set then any prime ideal Q of the ring $S_{n-1}(R)$ that is linked to P must be semistable (2) Any extended prime ideal of the ring $S_{n-1}(R)$ is link krull symmetric . Using induction hypothesis namely that the result is true for the case n-1, then by an argument similar to the case n=1, given above in corollary (6) (for this observe that we just have to use the above two properties only of the ring ring $S_{n-1}(R)$ ) we can prove the result for the case n . This completes the proof.